\documentclass[11pt]{article}%
\usepackage{amssymb}
\usepackage{graphicx}
\usepackage{amsmath}
\usepackage{amsthm}
\usepackage{amsfonts}%
\newtheorem{thm}{Theorem}[section]
\newtheorem{lemma}[thm]{Lemma}

\begin{document}

\title{Alternate Pacing of Border-Collision Period-Doubling Bifurcations}
\author{Xiaopeng Zhao$^{1}$\thanks{Corresponding author. Email: xzhao@duke.edu}~~and David G. Schaeffer$^{2}$\\Department of $^{1}$Biomedical Engineering and $^{2}$Mathematics \\and Center for Nonlinear and Complex Systems \\Duke University, Durham, North Carolina 27708}
\date{}
\maketitle

\begin{abstract}
Unlike classical bifurcations, border-collision bifurcations occur when, for
example, a fixed point of a continuous, piecewise $\mathcal{C}^{1}$ map
crosses a boundary in state space. Although classical bifurcations have been
much studied, border-collision bifurcations are not well understood. This
paper considers a particular class of border-collision bifurcations, i.e.,
border-collision period-doubling bifurcations. We apply a subharmonic
perturbation to the bifurcation parameter, which is also known as alternate
pacing, and we investigate the response under such pacing near the original
bifurcation point. The resulting behavior is characterized quantitatively by a
gain, which is the ratio of the response amplitude to the applied perturbation
amplitude. The gain in a border-collision period-doubling bifurcation has a
qualitatively different dependence on parameters from that of a classical
period-doubling bifurcation. Perhaps surprisingly, the differences are more
readily apparent if the gain is plotted vs. the perturbation amplitude (with
the bifurcation parameter fixed) than if plotted vs. the bifurcation parameter
(with the perturbation amplitude fixed). When this observation is exploited,
the gain under alternate pacing provides a useful experimental tool to
identify a border-collision period-doubling bifurcation.

\end{abstract}

\section{Introduction}

In contrast with the usual $\mathcal{C}^{1}$ context for classical
bifurcations \cite{Golubitsky}, the context of border-collision (B/C)
bifurcations is continuous, piecewise $\mathcal{C}^{1}$ maps
\cite{diBernardo06,BCbook}. The simplest normal form of a B/C bifurcation,
which we derive below, occurs for iteration of a piecewise linear map as
follows%
\begin{equation}
x_{n+1}=\left\{
\begin{array}
[c]{cc}%
A\,x_{n}+c\,\mu, & \text{if }x_{n}^{\left(  1\right)  }\geq0\\
B\,x_{n}+c\,\mu, & \text{if }x_{n}^{\left(  1\right)  }\leq0
\end{array}
\right.  , \label{eqn:normal_form}%
\end{equation}
where $x$ is an $m$-dimensional vector, $A$ and $B$ are $m\times m$ constant
matrices, $c$ is a $m$-dimensional constant vector, and the scalar $\mu$
represents a bifurcation parameter. Here, the term \emph{border} refers to
$x^{\left(  1\right)  }=0$, where the superscript indicates the first
component of $x$: i.e., the border is the surface across which derivatives of
the mapping jump.

Examples of B/C bifurcations have been observed in various engineering and
biological systems \cite{diBernardo06,nusse92,nusse94,Sun1995}. Different
phenomena in B/C bifurcations have been classified by a few authors
\cite{Banerjee1999,diBernardo1999,diBernardo06,Feigin1970,Feigin1974,Feigin1978,HassounehThesis}%
. Border-collision bifurcations may exhibit surprising phenomena such as
so-called instant chaos, which is a direct transition from a fixed point to a
chaotic attractor \cite{nusse94}. Alternatively, B/C bifurcations may bear
great similarity to their smooth counterparts, such as B/C period-doubling
bifurcations \cite{diBernardo06,BCbook}, which are the main focus of this paper.

\begin{figure}[ptb]
\centering\includegraphics[width=3.25in]{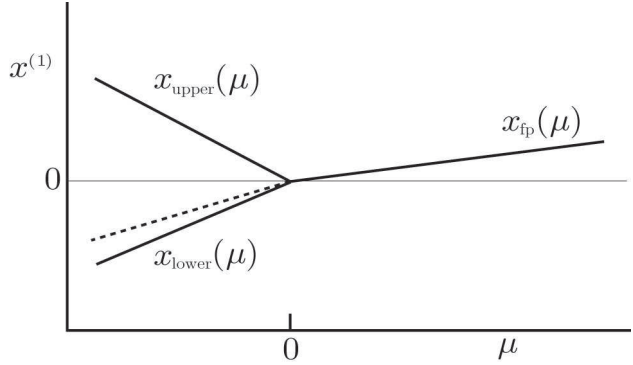}\caption{Schematic
bifurcation diagram of a border-collision period-doubling
bifurcation. Here, the thin horizontal line represents the
\emph{border}. Solid lines represent
stable solutions whereas the dashed line represents the unstable fixed point.}%
\label{fig:bif_diagram}%
\end{figure}

Unlike a classical period-doubling bifurcation that occurs when one eigenvalue
crosses the unit circle through $-1$, eigenvalues are not predictive for the
onset of a B/C period-doubling bifurcation. Instead, as indicated in Figure
\ref{fig:bif_diagram}, a B/C period-doubling bifurcation occurs when a branch
of fixed point \emph{collides} with a border. Despite such intrinsic
difference, it may be hard to distinguish between the two bifurcations,
especially based on necessarily limited data collected from experiments
\cite{Berger06}. This work investigates a possible technique to differentiate
between classical and B/C period-doubling bifurcations, inspired by a
phenomenon known as prebifurcation amplification previously observed in
classical period-doubling bifurcations \cite{vohra95,wiesenfeld86}. The term
\textquotedblleft prebifurcation amplification\textquotedblright\ refers to
the fact that, near the onset of a classical period-doubling bifurcation, a
subharmonic perturbation to the bifurcation parameter leads to amplified
disturbances in the system response. In the literature, such subharmonic
perturbations are also referred to as \emph{alternate pacing}
\cite{Berger06,Karma2004,zhao06prebif-smooth}. The amplification due to
alternate pacing of a classical bifurcation has been characterized by a gain
relating the disturbances in response to the perturbation amplitude
\cite{vohra95,wiesenfeld86,zhao06prebif-smooth}. In this paper, we compute the
gain for a B/C period-doubling bifurcation under alternate pacing.
Specifically, we find that the gain of a B/C bifurcation, as a function of the
bifurcation parameter and the perturbation amplitude, differs significantly
from that of a classical bifurcation.

In Section 2, we first investigate the conditions for the existence
and stability of a B/C period-doubling bifurcation. Then, the
response of a B/C period-doubling bifurcation to alternate pacing is
derived in Section 3. Section 4 quantifies the effect of
prebifurcation amplification through a gain. Finally, Section 5
provides a concluding discussion.

\section{Hypotheses and Notation}

\subsection{Normal form of a border-collision bifurcation}

We first derive the normal form (\ref{eqn:normal_form}). Consider iteration
under a general border-collision (B/C) map
\begin{equation}
z_{n+1}=f\left(  z_{n},\nu\right)  , \label{eqn:map_z}%
\end{equation}
where $z_{n}\in\mathbb{R}^{m}$, $\nu$ is a bifurcation parameter, and%
\begin{equation}
f\left(  z,\nu\right)  =\left\{
\begin{array}
[c]{cc}%
f_{\text{{\tiny A}}}\left(  z,\nu\right)  , & \text{if }\beta\left(
z,\nu\right)  \geq0\\
f_{\text{{\tiny B}}}\left(  z,\nu\right)  , & \text{if }\beta\left(
z,\nu\right)  \leq0
\end{array}
\right.  . \label{eqn:f_def}%
\end{equation}
To have a B/C bifurcation, we impose the following three conditions on the
above map.

\textbf{Condition 1:} Continuity on the boundary: $f_{\text{{\tiny A}}}\left(
z,\nu\right)  =f_{\text{{\tiny B}}}\left(  z,\nu\right)  $ whenever
$\beta\left(  z,\nu\right)  =0$.

\textbf{Condition 2:} For $\nu=\nu_{\text{bif}}$, (\ref{eqn:f_def}) has a
fixed point $z_{\text{bif}}$ on the boundary, i.e.,
\begin{equation}
\beta\left(  z_{\text{bif}},\nu_{\text{bif}}\right)  =0
\end{equation}
and $z_{\text{bif}}=f_{\text{{\tiny A}}}\left(  z_{\text{bif}},\nu
_{\text{bif}}\right)  =f_{\text{{\tiny B}}}\left(  z_{\text{bif}}%
,\nu_{\text{bif}}\right)  $.

\textbf{Condition 3:} The boundary is locally nonsingular; in symbols%
\begin{equation}
\nabla\beta\left\vert _{\text{bif}}\right.  \equiv\nabla\beta\left(
z_{\text{bif}},\nu_{\text{bif}}\right)  \neq0,
\end{equation}
where $\nabla$ indicates the vector of derivatives of $\beta$ with respective
to $z_{1},\ldots,z_{m};$ i.e., it does not include the $\nu$ derivative.

To derive the normal form, choose an invertible matrix $R$ s.t.
\begin{equation}
R\,\left(  \nabla\beta\left\vert _{\text{bif}}\right.  \right)  =e_{1},
\end{equation}
where $e_{1}=\left(  1,0,\ldots,0\right)  ^{T}$. Let $b$ be the vector%
\begin{equation}
b=\left(  \partial_{\nu}\beta\left\vert _{\text{bif}}\right.  \right)
\,e_{1}.
\end{equation}
Making a change of coordinates%
\begin{equation}
x=R\,\left(  z-z_{\text{bif}}\right)  +\left(  \nu-\nu_{\text{bif}}\right)
\,b,\,\,\,\,\,\,\,\,\mu=\nu-\nu_{\text{bif}} \label{eqn:newcoord_mu}%
\end{equation}
and dropping all higher order terms in (\ref{eqn:map_z}) yields a map that
assumes the form (\ref{eqn:normal_form}), where%
\begin{align}
A  &  =R\,\left(  Df_{\text{{\tiny A}}}\left\vert _{\text{bif}}\right.
\right)  \,R^{-1},\,\,\,\,\,\,\,\,B=R\,\left(  Df_{\text{{\tiny B}}}\left\vert
_{\text{bif}}\right.  \right)  \,R^{-1}\\
c  &  =b-A\,b-R\,\left(  \partial_{\nu}f_{\text{{\tiny A}}}\left\vert
_{\text{bif}}\right.  \right)  =b-B\,b-R\,\left(  \partial_{\nu}%
f_{\text{{\tiny B}}}\left\vert _{\text{bif}}\right.  \right)  .
\end{align}
Here, $Df$ denotes the differential of $f$, i.e., the $m\times m$ matrix of
partial derivative of $f$ with respective to $z_{1},\ldots z_{m}$. Note that
the two expressions of the vector $c$ are equal due to the continuity of $f$.
As will be seen later, condition 4 below requires that $c\neq0$.

\textbf{Remark:} The continuity condition (\#1 above) implies that $A\,x=B\,x$
if $x^{\left(  1\right)  }=0$. Thus, only nonzero column of $B-A$ is the first column.

\subsection{Conditions for the period-doubling case}

For the rest of the paper, we assume map (\ref{eqn:normal_form}) undergoes a
B/C period-doubling bifurcation as illustrated in Figure \ref{fig:bif_diagram}%
. Specifically, we assume that:

\textbf{Hypothesis 1:} (a) For $\mu>0$, map (\ref{eqn:normal_form}) has a
stable fixed point $x_{\text{fp}}\left(  \mu\right)  $ in the half space
$\left\{  x^{\left(  1\right)  }>0\right\}  $. (b) For $\mu<0$, there exists a
fixed point in the half space $\left\{  x^{\left(  1\right)  }<0\right\}  $
but it is unstable.

\textbf{Hypothesis 2:} For $\mu<0$, (\ref{eqn:normal_form}) has a stable
period-two orbit, with one point in $\left\{  x^{\left(  1\right)
}>0\right\}  $ and the other in $\left\{  x^{\left(  1\right)  }<0\right\}  $.

Conditions 4-6 below guarantee that (\ref{eqn:normal_form}) possesses this
structure. We do \emph{not} assume that the fixed point for $\mu>0$ is the
only attractor nor do we assume the period-two orbit for $\mu<0$ is the only
attractor. Indeed, the analysis of Banerjee and Grebogi \cite{Banerjee1999}
suggests that, more often than not, there will be additional attractors.

\subsubsection{{Analysis of fixed point solutions}}

When $\mu>0$, a fixed point of (\ref{eqn:normal_form}) in the half space
$\left\{  x^{\left(  1\right)  }>0\right\}  $ must satisfy
\begin{equation}
x_{\text{fp}}\left(  \mu\right)  =\mu X_{\text{fp}}, \label{eqn:x_fp}%
\end{equation}
where $X_{\text{fp}}$ is the constant vector%
\begin{equation}
X_{\text{fp}}=\left(  I-A\right)  ^{-1}\,c.
\end{equation}
By hypothesis 1 (a), we have the following condition.

\textbf{Condition 4:} The first component of the vector $X_{\text{fp}}$ is
positive, i.e., $X_{\text{fp}}^{\left(  1\right)  }>0$, and all eigenvalues of
$A$ are inside the unit circle, i.e., $\forall i,$ $\left\vert \lambda
_{i}\left(  A\right)  \right\vert <1$.

When $\mu<0$, the unstable fixed point in $\left\{  x^{\left(  1\right)
}<0\right\}  $ is given by
\begin{equation}
x_{\text{unstable}}\left(  \mu\right)  =\mu\,\left(  I-B\right)  ^{-1}\,c,
\end{equation}
and by hypothesis 1 (b) we have the following condition.

\textbf{Condition 5: }The first component of the vector\textbf{ }$\left(
I-B\right)  ^{-1}\,c$ is positive and at least one eigenvalue of $B$ satisfies
that $\left\vert \lambda_{i}\left(  B\right)  \right\vert >1$.

\subsubsection{{Analysis of period-two solutions}}

When $\mu<0$, let us write $x_{\text{upper}}\left(  \mu\right)  $,
$x_{\text{lower}}\left(  \mu\right)  $ for the stable period-two orbit of
(\ref{eqn:normal_form}), where as in Figure \ref{fig:bif_diagram}
$x_{\text{upper}}\left(  \mu\right)  $ is the point in $\left\{  x^{\left(
1\right)  }>0\right\}  $. These points satisfy%
\begin{equation}
x_{\text{upper}}\left(  \mu\right)  =B\,x_{\text{lower}}\left(  \mu\right)
+\mu\,c,\,\,\,\,\,\,\,x_{\text{lower}}\left(  \mu\right)  =A\,x_{\text{upper}%
}\left(  \mu\right)  +\mu\,c \label{eqn:xul_eqn}%
\end{equation}
It follows from the above equations that
\begin{equation}
x_{\text{upper}}=\mu\,X_{\text{upper}},\,\,\,\,\,\,\,x_{\text{lower}}%
=\mu\,X_{\text{lower}}, \label{eqn:x_u&l}%
\end{equation}
where%
\begin{equation}
X_{\text{upper}}=\left(  I-B\,A\right)  ^{-1}\,\left(  I+B\right)
\,c,\,\,\,\,\,\,\,X_{\text{lower}}=\left(  I-A\,B\right)  ^{-1}\,\left(
I+A\right)  \,c.
\end{equation}
For these calculations to be consistent and for the orbit to be stable, we
need the following condition.

\textbf{Condition 6:} $X_{\text{upper}}^{\left(  1\right)  }<0$,
$X_{\text{lower}}^{\left(  1\right)  }>0$, and $\forall i,$ $\left\vert
\lambda_{i}\left(  A\,B\right)  \right\vert <1$. (Note that matrices $A\,B$
and $B\,A$ have same characteristic polynomials and hence have the same eigenvalues.)

\section{The response to alternate pacing}

Alternate pacing imposes a subharmonic perturbation to the bifurcation
parameter $\mu$, rendering the map (\ref{eqn:normal_form}) as%
\begin{equation}
x_{n+1}=\left\{
\begin{array}
[c]{cc}%
A\,x_{n}+\left[  \mu+\left(  -1\right)  ^{n}\,\delta\right]  \,c, & \text{if
}x_{n}^{\left(  1\right)  }>0\\
B\,x_{n}+\left[  \mu+\left(  -1\right)  ^{n}\,\delta\right]  \,c, & \text{if
}x_{n}^{\left(  1\right)  }<0
\end{array}
\right.  , \label{eqn:app}%
\end{equation}
For use below, we define%
\begin{equation}
d=\left(  I+A\right)  ^{-1}\,c; \label{eqn:d}%
\end{equation}
to avoid a degenerate response to alternate pacing, we impose the following condition.

\textbf{Condition 7:} The first component of the vector $d$ is nonzero, i.e.,
$d^{\left(  1\right)  }\neq0$.

\begin{figure}[ptb]
\centering\includegraphics[width=4in]{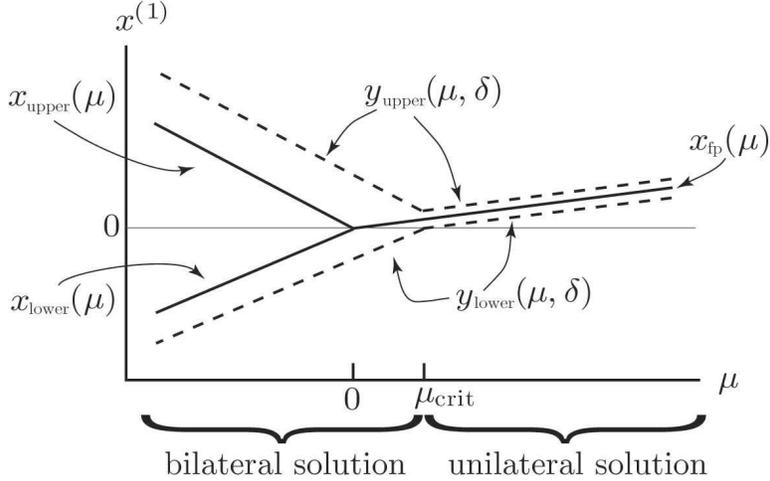}\caption{Unilateral
and bilateral solutions of alternate pacing (dashed) compared to the
solutions
when $\delta=0$ (solid). The thin thin horizontal line represents the border.}%
\label{fig:app_bifdiagram}%
\end{figure}

We are interested in period-two solutions of (\ref{eqn:app}), whose two points
are denoted by $y_{\text{upper}}\left(  \mu,\delta\right)  $ and
$y_{\text{lower}}\left(  \mu,\delta\right)  $. As depicted in Figure
\ref{fig:app_bifdiagram} and will be derived in the following, the character
of the period-two solution depends on $\mu$. 1.) For $\mu>\mu_{\text{crit}%
}\left(  \delta\right)  $, where $\mu_{\text{crit}}\left(  \delta\right)  $ is
defined in (\ref{eqn:mu_crit}) below, the response is a unilateral solution:
i.e., both points are above the border, or in symbols $y_{\text{upper}%
}^{\left(  1\right)  }>y_{\text{lower}}^{\left(  1\right)  }>0$; and
2.) For $\mu<\mu_{\text{crit}}\left(  \delta\right)  $, the response
is a bilateral solution: i.e., one point is above and the other is
below the border, or in symbols $y_{\text{upper}}^{\left(  1\right)
}>0>y_{\text{lower}}^{\left( 1\right)  }$.

\subsection{The unilateral solution}

To compute the unilateral solution, we first need to determine whether
$y_{\text{upper}}\left(  \mu,\delta\right)  $ occurs in the iteration for $n$
even or odd. To this end, we temporarily denote the two points as
$y_{\text{odd}}\left(  \mu,\delta\right)  $ and $y_{\text{even}}\left(
\mu,\delta\right)  $. In the unilateral case, i.e. $y_{\text{odd}}^{\left(
1\right)  }$, $y_{\text{even}}^{\left(  1\right)  }>0$, equation
(\ref{eqn:app})\ leads to%
\begin{equation}
y_{\text{odd}}=A\,y_{\text{even}}+\left(  \mu+\delta\right)
c,\,\,\,\,\,\,\,y_{\text{even}}=A\,y_{\text{odd}}+\left(  \mu-\delta\right)
c. \label{eqn:yeven&odd_eqn}%
\end{equation}
Solving the above equations yields%
\begin{equation}
y_{\text{odd}}\left(  \mu,\delta\right)  =x_{\text{fp}}\left(  \mu\right)
+\delta\,d,\,\,\,\,\,\,\,y_{\text{even}}\left(  \mu,\delta\right)
=x_{\text{fp}}\left(  \mu\right)  -\delta\,d, \label{eqn:yeven&odd}%
\end{equation}
where $x_{\text{fp}}$ and $d$ are computed in (\ref{eqn:x_fp}) and
(\ref{eqn:d}), respectively. By adopting the sign convention%
\begin{equation}
\text{sign}~\delta=\text{sign}~d^{\left(  1\right)  }, \label{eqn:sign}%
\end{equation}
we can arrange that $y_{\text{upper}}$ and $y_{\text{lower}}$ correspond to
odd- and even-numbered beats, respectively; that is%
\begin{equation}
y_{\text{upper}}\left(  \mu,\delta\right)  =x_{\text{fp}}\left(  \mu\right)
+\delta\,d,\,\,\,\,\,\,\,y_{\text{lower}}\left(  \mu,\delta\right)
=x_{\text{fp}}\left(  \mu\right)  -\delta\,d. \label{eqn:uni_sol}%
\end{equation}
In the remainder of the paper, we retain the notation $y_{\text{upper}}$ and
$y_{\text{lower}}$, assuming the sign convention (\ref{eqn:sign}).

The unilateral solution is valid provided $y_{\text{lower}}^{\left(  1\right)
}\geq0$, that is
\begin{equation}
\mu\geq\rho\,\left\vert \delta\right\vert \equiv\mu_{\text{crit}}\left(
\delta\right)  ,\label{eqn:mu_crit}%
\end{equation}
where%
\begin{equation}
\rho=\frac{\left\vert d^{\left(  1\right)  }\right\vert }{X_{\text{fp}%
}^{\left(  1\right)  }}>0\text{.}\label{eqn:rho}%
\end{equation}
At the critical value of $\mu$, we have%
\begin{align}
y_{\text{upper}}\left(  \mu_{\text{crit}}\left(  \delta\right)  ,\delta
\right)   &  =x_{\text{fp}}\left(  \mu_{\text{crit}}\left(  \delta\right)
\right)  +\delta\,d\equiv y_{\text{u-crit}}\left(  \delta\right)
\label{eqn:yu_crit}\\
y_{\text{lower}}\left(  \mu_{\text{crit}}\left(  \delta\right)  ,\delta
\right)   &  =x_{\text{fp}}\left(  \mu_{\text{crit}}\left(  \delta\right)
\right)  -\delta\,d\equiv y_{\text{l-crit}}\left(  \delta\right)
.\label{eqn:yl_crit}%
\end{align}
Note that $y_{\text{u-crit}}^{\left(  1\right)  }>0$ and $y_{\text{l-crit}%
}^{\left(  1\right)  }=0$.

\subsection{The bilateral solution}

We now compute the bilateral solution that evolves continuously from the
unilateral solution\footnote{In the appendix, we show that for $\mu$
sufficiently large negative there is also a second bilateral solution.} as
$\mu$ crosses $\mu_{\text{crit}}\left(  \delta\right)  $. Thus, for $\mu
<\mu_{\text{crit}}\left(  \delta\right)  $, we look for a period-two orbit of
(\ref{eqn:app}) such that
\begin{equation}
y_{\text{upper}}^{\left(  1\right)  }\left(  \mu,\delta\right)
>0>y_{\text{lower}}^{\left(  1\right)  }\left(  \mu,\delta\right)  ,
\label{eqn:bilateral}%
\end{equation}
where $y_{\text{upper}}$ occurs for $n$ odd. By (\ref{eqn:app}), these satisfy%

\begin{equation}
y_{\text{upper}}=B\,y_{\text{lower}}+\left(  \mu+\delta\right)
c,\,\,\,\,\,\,\,y_{\text{lower}}=A\,y_{\text{upper}}+\left(  \mu
-\delta\right)  c. \label{eqn:bilateral_eqn}%
\end{equation}
Comparing these equations with (\ref{eqn:xul_eqn}), we see that the solution
of (\ref{eqn:bilateral_eqn}) may be written as%

\begin{equation}
y_{\text{upper}}\left(  \mu,\delta\right)  =x_{\text{upper}}\left(
\mu\right)  +\delta\,s_{\text{upper}},\,\,\,\,\,\,\,y_{\text{lower}}\left(
\mu,\delta\right)  =x_{\text{lower}}\left(  \mu\right)  -\delta
\,s_{\text{lower}}, \label{eqn:bilateral_sol}%
\end{equation}
where $x_{\text{upper}}$ and $x_{\text{lower}}$ have been computed in
(\ref{eqn:x_u&l}), and the $\mu$-independent shifts $s_{\text{upper}}$,
$s_{\text{lower}}$ are given by%
\begin{equation}
s_{\text{upper}}=\left(  I-B\,A\right)  ^{-1}\,\left(  I-B\right)
\,c,\,\,\,\,\,\,\,s_{\text{lower}}=\left(  I-A\,B\right)  ^{-1}\,\left(
I-A\right)  \,c. \label{eqn:s_u&l}%
\end{equation}

We now show that this solution is consistent with the assumption
(\ref{eqn:bilateral}). We claim that when $\mu=\mu_{\text{crit}}\left(
\delta\right)  $, the solution (\ref{eqn:yu_crit},\ref{eqn:yl_crit}) of
(\ref{eqn:yeven&odd_eqn}) for this value of $\mu$ also satisfies
(\ref{eqn:bilateral_eqn}); this claim follows from the observation that, since
$y_{\text{l-crit}}^{\left(  1\right)  }=0$, we have $A\,y_{\text{l-crit}%
}=B\,y_{\text{l-crit}}$. Using the fact that $x_{\text{upper}}$ and
$x_{\text{lower}}$ are linear in $\mu$, we may rewrite
(\ref{eqn:bilateral_sol}) as
\begin{align}
y_{\text{upper}}\left(  \mu,\delta\right)   &  =y_{\text{upper}}\left(
\mu_{\text{crit}}\left(  \delta\right)  ,\delta\right)  +x_{\text{upper}%
}\left(  \mu\right)  -x_{\text{upper}}\left(  \mu_{\text{crit}}\left(
\delta\right)  \right) \nonumber\\
&  =y_{\text{u-crit}}\left(  \delta\right)  +\left(  \mu-\mu_{\text{crit}%
}\left(  \delta\right)  \right)  \,X_{\text{upper}} \label{eqn:yupper_re}%
\end{align}
and%
\begin{equation}
y_{\text{lower}}\left(  \mu,\delta\right)  =y_{\text{l-crit}}\left(
\delta\right)  +\left(  \mu-\mu_{\text{crit}}\left(  \delta\right)  \right)
\,X_{\text{lower}}. \label{eqn:ylower_re}%
\end{equation}
It now follows from Condition 6 that (\ref{eqn:bilateral}) holds for all
$\mu<\mu_{\text{crit}}\left(  \delta\right)  $.

\subsection{Stability of perturbed iterations}

The unperturbed map (\ref{eqn:normal_form}), without alternate pacing as in
(\ref{eqn:app}), may suffer a so-called dangerous bifurcation
\cite{HassounehPRL} at $\mu=0$. To elaborate: suppose that at some positive
value of $\mu$ the (unperturbed) iteration is locked onto $x_{\text{fp}%
}\left(  \mu\right)  $. If $\mu$ is decreased quasistatically but remains
positive, the solution will follow the stable fixed point $x_{\text{fp}%
}\left(  \mu\right)  $. However, as $\mu\rightarrow0$, the basin of attraction
of the fixed point might shrink to nil. In this case, the iteration becomes
indeterminant as $\mu$ crosses $0$, i.e., as a result of infinitesimal
perturbations, the solution might shift to a different attractor for $\mu<0$
rather than to the period-two solution.

One may hope that iterations of (\ref{eqn:app}), i.e., the map under alternate
pacing, are better behaved near $\mu=\mu_{\text{crit}}\left(  \delta\right)  $
when they cross the boundary plane $\left\{  x^{\left(  1\right)  }=0\right\}
$. However, we are unable to guarantee this without imposing an additional
hypothesis, as follows.

\textbf{Condition 8:} There exists an invertible matrix $S$ s.t. $\left\Vert
S\,A^{2}\,S^{-1}\right\Vert <1$ and $\left\Vert S\,A\,B\,S^{-1}\right\Vert
<1$. In other words, the eigenvectors for $A^{2}$ and $A\,B$ are not too
different. For example, this condition is satisfied if both $A^{2}$ and $A\,B$
are both diagonalizable using $S$.

In the following lemma, we consider the composition of (\ref{eqn:app}) with
itself when $\mu\approx\mu_{\text{crit}}\left(  \delta\right)  $, starting
with $n$ even from an $x$ near the even iterates $y_{\text{lower}}$ of the
period-two orbit, for which we have $y_{\text{lower}}^{\left(  1\right)
}\left(  \mu,\delta\right)  \approx y_{\text{l-crit}}^{\left(  1\right)  }=0$.
Thus, both options in (\ref{eqn:app}) must be considered for the first
application of (\ref{eqn:app}). However, the image of $x$ will be close to the
odd iterates, for which we have $y_{\text{upper}}^{\left(  1\right)  }\left(
\mu,\delta\right)  \approx y_{\text{u-crit}}^{\left(  1\right)  }>0$. Thus,
provided $\mu$ is sufficiently close to $\mu_{\text{crit}}\left(
\delta\right)  $ and $x$ is sufficiently close to $y_{\text{lower}}^{\left(
1\right)  }\left(  \mu,\delta\right)  $, only the upper option in
(\ref{eqn:app}) will occurs in the second application of (\ref{eqn:app}).
Thus, if $\mu$ and $x$ are so restricted, the composition of (\ref{eqn:app})
with itself is given by
\begin{equation}
F\left(  x\right)  =\left\{
\begin{array}
[c]{cc}%
A\left[  A\,x+\left(  \mu+\delta\right)  c\right]  +\left(  \mu-\delta\right)
c & \text{if }x^{\left(  1\right)  }\geq0\\
A\left[  B\,x+\left(  \mu+\delta\right)  c\right]  +\left(  \mu-\delta\right)
c & \text{if }x^{\left(  1\right)  }\leq0
\end{array}
\right.  . \label{eqn:F}%
\end{equation}

\begin{lemma}
$F$ maps a neighborhood of $y_{\text{lower}}\left(  \mu,\delta\right)  $ into
itself and, with respect to an appropriate norm, is a contraction there.
\end{lemma}

\textbf{Proof:} Given a nonsingular matrix $S$, we define the $S$-norm of a
vector $v$ as
\begin{equation}
\left\Vert v\right\Vert _{S}=\left\Vert S\,v\right\Vert
\end{equation}
and that of a matrix $M$ as%
\begin{equation}
\left\Vert M\right\Vert _{S}=\left\Vert S\,M\,S^{-1}\right\Vert .
\end{equation}
For the matrix $S$ in Condition 8, we have $\theta=\max\left\{  \left\Vert
A^{2}\right\Vert _{S},\left\Vert A\,B\right\Vert _{S}\right\}  <1$.

To prove $F$ in (\ref{eqn:F}) is a contraction, we consider $x$ and $y$ in a
neighborhood of $y_{\text{lower}}$ and treat the following three cases separately.

1.) $x^{\left(  1\right)  }$, $y^{\left(  1\right)  }\geq0$. It follows from
(\ref{eqn:F}) that%
\begin{equation}
\left\Vert F\left(  x\right)  -F\left(  y\right)  \right\Vert _{S}=\left\Vert
A^{2}\left(  x-y\right)  \right\Vert _{S}\leq\left\Vert A^{2}\right\Vert
_{S}\left\Vert \left(  x-y\right)  \right\Vert _{S}\leq\theta\left\Vert
x-y\right\Vert _{S}%
\end{equation}

2.) $x^{\left(  1\right)  }$, $y^{\left(  1\right)  }\leq0$. It follows from
(\ref{eqn:F}) that%
\begin{equation}
\left\Vert F\left(  x\right)  -F\left(  y\right)  \right\Vert _{S}=\left\Vert
A\,B\,\left(  x-y\right)  \right\Vert _{S}\leq\left\Vert A\,B\right\Vert
_{S}\left\Vert \left(  x-y\right)  \right\Vert _{S}\leq\theta\left\Vert
x-y\right\Vert _{S}%
\end{equation}

3.) $x^{\left(  1\right)  }\geq0$ and $y^{\left(  1\right)  }\leq0$. Let $z$
be the point where the line between $x$ and $y$ intersects the plane $\left\{
z^{\left(  1\right)  }=0\right\}  $. This choice of $z$ implies that
$A\,z=B\,z$ and $\left\Vert x-z\right\Vert _{S}+\left\Vert z-y\right\Vert
_{S}=\left\Vert x-y\right\Vert _{S}$. It follows from (\ref{eqn:F}) that
\begin{align}
\left\Vert F\left(  x\right)  -F\left(  y\right)  \right\Vert _{S} &
=\left\Vert A^{2}\,x-A\,B\,y\right\Vert _{S}=\left\Vert A^{2}\left(
x-z\right)  +A\,B\,\left(  z-y\right)  \right\Vert _{S}\\
&  \leq\left\Vert A^{2}\right\Vert _{S}\left\Vert x-z\right\Vert
_{S}+\left\Vert A\,B\right\Vert _{S}\left\Vert y-z\right\Vert _{S}\\
&  \leq\theta\left\Vert x-z\right\Vert _{S}+\theta\left\Vert y-z\right\Vert
_{S}=\theta\left\Vert x-y\right\Vert _{S}.
\end{align}
Moreover, since $y_{\text{lower}}$ is a fixed point, $F$ maps a neighborhood
of $y_{\text{lower}}$ into itself.\vrule height 6pt width 6pt depth 1pt

\begin{thm}
Condition 8 implies that the period-two orbit $y_{\text{upper}}$,
$y_{\text{lower}}$ is stable.
\end{thm}

\textbf{Proof:} For $\mu>\mu_{\text{crit}}$ and $\mu<\mu_{\text{crit}}$,
stability of the iterations is guaranteed because $\forall i$, $\left\vert
\lambda_{i}\left(  A^{2}\right)  \right\vert <1$ and $\left\vert \lambda
_{i}\left(  B\,A\right)  \right\vert <1$, respectively. For $\mu\approx
\mu_{\text{crit}}$, stability is guaranteed because $F$ is a contraction, as
shown in the proof of Lemma 3.1. \vrule height 6pt width 6pt depth 1pt

\section{Prebifurcation gain}

\subsection{Analysis}

Near the onset of a classical period-doubling bifurcation, disturbances in the
bifurcation parameter may cause amplified disturbances in the system response,
a phenomenon known as prebifurcation amplification. Various authors studied
prebifurcation amplification for a classical period-doubling bifurcation,
using a gain, which is the ratio of the response magnitude to the applied
perturbation amplitude \cite{surovyatkina2005,vohra95,wiesenfeld86}. An
analytical formula of the gain, including crucial higher-order terms, was
derived in \cite{zhao06prebif-smooth}. It is interesting to compute the gain
for a B/C period-doubling bifurcation and compare with that of a classical
period-doubling bifurcation. To this end, let us define a gain using the 1st
component of the solution%
\begin{equation}
\Gamma\equiv\frac{y_{\text{upper}}^{\left(  1\right)  }-y_{\text{lower}%
}^{\left(  1\right)  }}{2\,\left\vert \delta\right\vert }.
\label{eqn:Gamma_def}%
\end{equation}
Then, the prebifurcation gain (when $\mu\geq0$) is described by the following theorem.

\begin{thm}
When $\mu\geq0$ and $\delta\neq0$, the gain satisfies%
\begin{equation}
\Gamma=\left\{
\begin{array}
[c]{cc}%
\Gamma_{\text{const}}, & \text{if }\mu\geq\rho\left\vert \delta\right\vert \\
\Gamma_{\text{const}}+\gamma\left(  \rho-\frac{\mu}{\left\vert \delta
\right\vert }\right)   & \text{if }\mu\leq\rho\left\vert \delta\right\vert
\end{array}
\right.  ,\label{eqn:theory}%
\end{equation}
where $\Gamma_{\text{const}}$ and $\gamma$ are positive constant, and $\rho$
is given by (\ref{eqn:rho}).
\end{thm}

\begin{figure}[ptb]
\centering%
\begin{tabular}
[c]{c}%
\includegraphics[width=4.75in]{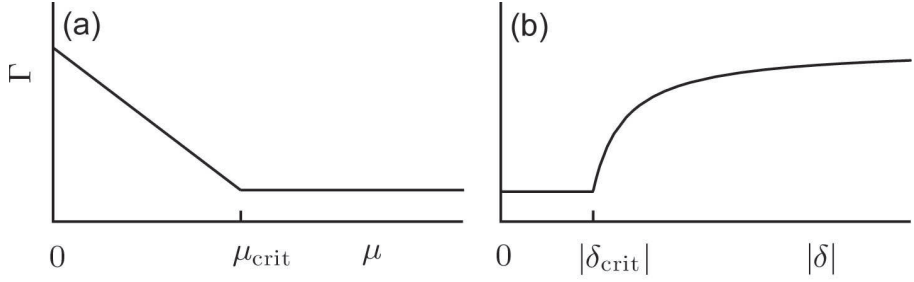}
\end{tabular}
\caption{Schematic of the variation of the canonical gain (\ref{eqn:Gamma_def}%
) under changes in parameters: (a)\ $\Gamma$ vs. $\mu$ when $\delta$ stays
constant and (b) $\Gamma$ vs. $\delta$ when $\mu$ stays constant.}%
\label{fig:gain_1}%
\end{figure}

\textbf{Proof: }When $\mu\geq\rho\left\vert \delta\right\vert $, the response
to alternate pacing is an unilateral solution, as shown in the previous
section. It follows from Equation (\ref{eqn:uni_sol}) that
\begin{equation}
\Delta y=y_{\text{upper}}-y_{\text{lower}}=2\,\delta
\,d.\label{eqn:type1_deltay}%
\end{equation}
By the definition (\ref{eqn:Gamma_def}) and recalling the sign convention
(\ref{eqn:sign}), the gain is
\begin{equation}
\Gamma=\left\vert d^{\left(  1\right)  }\right\vert \equiv\Gamma
_{\text{const}}.\label{eqn:Gamma_const}%
\end{equation}
Therefore, if $\mu$ is varied and $\delta$ is held fixed, the gain stays
constant as long as $\mu\geq\mu_{\text{crit}}\left(  \delta\right)
=\rho\left\vert \delta\right\vert $, see Figure \ref{fig:gain_1} (a).
Similarly, if $\delta$ is varied and $\mu$ is held fixed, the gain stays
constant as long as $\left\vert \delta\right\vert \leq\delta_{\text{crit}%
}\left(  \mu\right)  $, where $\delta_{\text{crit}}\left(  \mu\right)
=\mu/\rho$, see Figure \ref{fig:gain_1} (b).

When $\mu\leq\rho\left\vert \delta\right\vert $, the response to alternate
pacing is a bilateral solution, as shown in the previous section. It follows
from Equations (\ref{eqn:yupper_re}) and (\ref{eqn:ylower_re}) that
\begin{align}
\Delta y &  =y_{\text{upper}}-y_{\text{lower}}\nonumber\\
&  =y_{\text{u-crit}}-y_{\text{l-crit}}+\left(  \mu_{\text{crit}}\left(
\delta\right)  -\mu\right)  \,\left(  X_{\text{upper}}-X_{\text{lower}%
}\right)  \label{eqn:type2_deltay}\\
&  =2\,\delta\,d+\left(  \mu_{\text{crit}}\left(  \delta\right)  -\mu\right)
\,\left(  X_{\text{upper}}-X_{\text{lower}}\right)  ,\nonumber
\end{align}
where we have used (\ref{eqn:yu_crit}) and (\ref{eqn:yl_crit}) for
$y_{\text{u-crit}}$ and $y_{\text{l-crit}}$. Recalling the definition
(\ref{eqn:Gamma_def}) and using (\ref{eqn:Gamma_const}), the gain is
\begin{equation}
\Gamma=\Gamma_{\text{const}}+\gamma\left(  \rho-\frac{\mu}{\left\vert
\delta\right\vert }\right)  ,
\end{equation}
where%
\begin{equation}
\gamma=\frac{1}{2}\left(  X_{\text{upper}}^{\left(  1\right)  }%
-X_{\text{lower}}^{\left(  1\right)  }\right)  \text{.}%
\end{equation}
Here, $\gamma>0$, since $X_{\text{upper}}^{\left(  1\right)  }>0$ and
$X_{\text{lower}}^{\left(  1\right)  }<0$. Therefore, for a constant $\delta$,
the gain increases as $\mu$ decreases from $\mu_{\text{crit}}=\rho\left\vert
\delta\right\vert $, see Figure \ref{fig:gain_1} (a). On the other hand, for a
constant $\mu$, the gain increases as $\left\vert \delta\right\vert $
increases from $\left\vert \delta_{\text{crit}}\right\vert =\mu/\rho$, see
Figure \ref{fig:gain_1} (b).\vrule height 6pt width 6pt depth 1pt

The precise relations (\ref{eqn:theory}), which hold for all $\mu$ and
$\delta$, result from the normal form (\ref{eqn:normal_form}), in which
higher-order terms have been dropped. If such terms are present, then
(\ref{eqn:theory}) will describe the qualitative behavior of the gain near the
bifurcation point for small perturbation amplitude. In Section 4.2 we present
a numerical example with the B/C period-doubling bifurcation in the cardiac
model of Sun et al. \cite{Sun1995}; this example shows that the leading-order
analysis accurately capture the behavior of the model.

\textbf{Remark: }Analogous to $\Gamma$, a generalized gain can be defined
using an arbitrary component of $y$ as
\begin{equation}
g^{\left(  i\right)  }\equiv\frac{y_{\text{upper}}^{\left(  i\right)
}-y_{\text{lower}}^{\left(  i\right)  }}{2\,\left\vert \delta\right\vert }.
\end{equation}
It follows from (\ref{eqn:type1_deltay}) and (\ref{eqn:type2_deltay}) that
\begin{equation}
g^{\left(  i\right)  }=\left\{
\begin{array}
[c]{cc}%
g_{\text{const}}^{\left(  i\right)  }, & \text{if }\mu\geq\rho\left\vert
\delta\right\vert \\
g_{\text{const}}^{\left(  i\right)  }+\left(  \rho-\frac{\mu}{\left\vert
\delta\right\vert }\right)  \,k^{\left(  i\right)  } & \text{if }\mu\leq
\rho\left\vert \delta\right\vert
\end{array}
\right.  ,
\end{equation}
where
\begin{equation}
g_{\text{const}}^{\left(  i\right)  }=\left(  \text{sign\thinspace}d^{\left(
1\right)  }\right)  \,\,d^{\left(  i\right)  },\,\,\,\,\,\,\,\,k^{\left(
i\right)  }=\frac{1}{2}\left(  X_{\text{upper}}^{\left(  i\right)
}-X_{\text{lower}}^{\left(  i\right)  }\right)  .
\end{equation}
As we saw above, $g_{\text{const}}^{\left(  1\right)  }=\Gamma_{\text{const}%
}>0$ and $k^{\left(  1\right)  }=\gamma>0$. However, if $i\neq1$, the signs of
$g_{\text{const}}^{\left(  i\right)  }$ and $k^{\left(  i\right)  }$ are not
uniquely determined. Thus, the generalized gain may have either sign and, even
if positive, may exhibit monotonicity opposite to that in Figure 3. However,
all gains are constant when $\mu\geq\rho\left\vert \delta\right\vert $.

\subsection{A numerical example}

\begin{figure}[ptb]
\centering
\begin{tabular}
[c]{c}%
\includegraphics[width=5.5in]{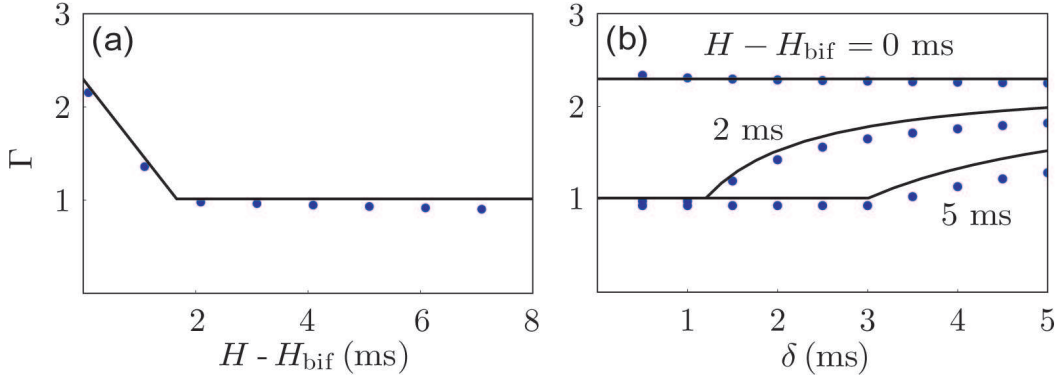}
\end{tabular}
\caption{Alternate pacing of the model of Sun \emph{et al.}
(\ref{eqn:sunmodel}): (a) variation of the gain $\Gamma$ under changes of $H$
when $\delta=1$ ms and (b) variation of the gain $\Gamma$ under changes of
$\delta$ when $H-H_{\text{bif}}=0$, $2$, and $5$ ms. Here, the solid curves
indicate theoretical results and the dots indicate numerical simulations.}%
\label{fig:sun}%
\end{figure}

Sun \emph{et al.} \cite{Sun1995} presented a theoretical model for
atrioventricular nodal conduction, which accurately predicts a variety of
experimentally observed cardiac rhythms. It was shown in \cite{HassounehIJBC}
that the model of Sun \emph{et al.} exhibits a border-collision
period-doubling bifurcation. Following the notation of \cite{HassounehIJBC},
the model iterates the atrial His interval, $A$, and the drift in the nodal
conduction time, $R$, as follows
\begin{equation}
\left(
\begin{array}
[c]{c}%
A_{n+1}\\
R_{n+1}%
\end{array}
\right)  =\left\{
\begin{array}
[c]{cc}%
\left(
\begin{array}
[c]{c}%
A_{\min}+R_{n+1}+\left(  201-0.7A_{n}\right)  e^{-H/\tau_{\text{rec}}}\\
R_{n}e^{-\left(  A_{n}+H\right)  /\tau_{\text{fat}}}+\gamma e^{-H/\tau
_{\text{fat}}}%
\end{array}
\right)  , & \text{if }A_{n}\leq130\\
& \\
\left(
\begin{array}
[c]{c}%
A_{\min}+R_{n+1}+\left(  500-3.0A_{n}\right)  e^{-H/\tau_{\text{rec}}}\\
R_{n}e^{-\left(  A_{n}+H\right)  /\tau_{\text{fat}}}+\gamma e^{-H/\tau
_{\text{fat}}}%
\end{array}
\right)  , & \text{if }A_{n}\geq130
\end{array}
\right.  , \label{eqn:sunmodel}%
\end{equation}
where $A_{\text{min}}=33$ ms, $\tau_{\text{rec}}=70$ ms, $\tau_{\text{fat}%
}=30\,000$ ms, and $\gamma=0.3$ ms. The bifurcation parameter, $H$, represents
the time interval between bundle of His activation and the subsequent
activation. Under variation of $H$, a border-collision period-doubling
bifurcation occurs at $H_{\text{bif}}=56.9078$ ms, where of course $A=130$ at
the bifurcation point. We apply alternate pacing to the above model by
perturbing $H$ with $\left(  -1\right)  ^{k}\delta$. Figure \ref{fig:sun}
shows the relation between the gain $\Gamma$ and the parameters $H$ and
$\delta$, computed from our theory as well as from numerical simulations. The
agreement between the theoretical and the numerical results is good.

\section{Summary and discussion}

A border-collision (B/C) bifurcation occurs when the fixed point of an
iterated map encounters a border on which the underlying map is continuous but
its derivatives jump. Using the normal form (\ref{eqn:normal_form}), we
presented the conditions for a B/C period-doubling bifurcation in multiple
dimensional maps. Under these conditions, we computed the response of the map
under a subharmonic perturbation to the bifurcation parameter, i.e. alternate
pacing. The response under alternate pacing can be either unilateral or
bilateral. As the names suggest, the two points of a unilateral solution lie
on the same side of the border, while those of a bilateral solution lie on
different sides of the border. The response to alternate pacing is quantified
through the definition of a gain. We found that the gain is a piecewise smooth
function of the bifurcation parameter ($\mu$) and the perturbation amplitude
($\left\vert \delta\right\vert $), with qualitatively different behaviors in
different parameter regions.

\begin{figure}[ptb]
\centering%
\begin{tabular}
[c]{c}%
\includegraphics[width=4.75in]{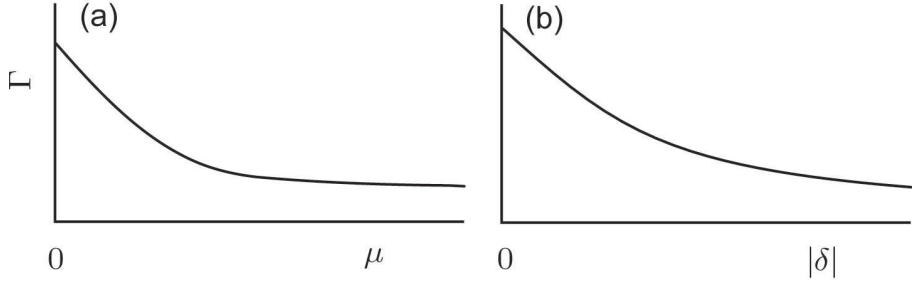}
\end{tabular}
\caption{Schematic of the variation of the gain of a classical bifurcation
under changes in parameters: (a)\ $\Gamma$ vs. $\mu$ when $\delta$ stays
constant and (b) $\Gamma$ vs. $\delta$ when $\mu$ stays constant.}%
\label{fig:gain_smooth}%
\end{figure}

Most importantly, the gain in a border-collision bifurcation differs
qualitatively from that in a classical bifurcation. For example, it follows
from \cite{zhao06prebif-smooth} that, after appropriate rescaling of the
parameters, the gain of a classical period-doubling bifurcation satisfies the
relation,%
\begin{equation}
\delta^{2}\,\Gamma^{3}+\mu\,\Gamma-1=0. \label{eqn:gain_smooth}%
\end{equation}
Figure \ref{fig:gain_smooth} (a) and (b) schematically show the dependence of
$\Gamma$ on $\mu$ and $\Gamma$ on $\delta$, respectively. Comparison of (\ref{eqn:theory}) and (\ref{eqn:gain_smooth})
reveals two differences between the gain in a classical bifurcation
and that in a border-collision bifurcation. (1) The gain of a
classical bifurcation tends to infinity\footnote{Moreover, the rate
of divergence as the parameters $\left(  \mu,\delta\right)  $ tend
to $\left(  0,0\right)  $ depends on the path taken. For example,
when $\delta$ is extremely small, the gain tends to infinity as
$\mu^{-1}$; on the other hand, when $\mu=0$, the gain tends to
infinity as $\delta^{-2/3}$.} as $\left(  \mu,\delta\right)  $
approaches $\left(  0,0\right)  $; by contrast, the gain of a
border-collision bifurcation is bound for all $\delta\neq0$. (2) The
gain of a classical bifurcation varies smoothly under changes in
system parameters while that of a border-collision bifurcation
undergoes a nonsmooth variation as parameters cross a boundary in
the parameter space.

Although the two cases exhibit these differences, if $\Gamma$ is plotted as a
function of $\mu$ (with $\delta$ held fixed), it may be difficult to
distinguish them based on data from experiments, because only discrete points
can be sampled and these are subject to experimental errors (cf. Figures
\ref{fig:gain_1} (a) and \ref{fig:gain_smooth} (a)). However, if $\Gamma$ is
plotted as a function of $\delta$ (with $\mu$ held fixed), the distinction
between the two bifurcation types is evident even with just a few data points
and in the presence of experimental noise (cf. Figures \ref{fig:gain_1} (b)
and \ref{fig:gain_smooth} (b)). Indeed, the theory in this paper guided the
experiments \cite{Berger06} in identifying the type of a bifurcation in a
paced cardiac tissue.

\begin{center}
\textbf{Acknowledgments}
\end{center}

Support of the National Institutes of Health under grant 1R01-HL-72831 and the
National Science Foundation under grants DMS-9983320 and PHY-0243584 is
gratefully acknowledged. The authors are also grateful to Carolyn Berger,
Daniel J. Gauthier, and Wanda Krassowska for their insightful discussion.

\bigskip\appendix{}

\section*{Appendix: The out-of-phase bilateral solution}

Recall that, for the unilateral solution under alternate pacing,
$y_{\text{upper}}$ occurs for $n$ odd due to the sign convention
(\ref{eqn:sign}), and this behavior continues for the bilateral solution
computed in Section 3.2; i.e. this bilateral solution is in phase with the
unilateral solution. Thus, we refer to such bilateral solution as an
\emph{in-phase} bilateral solution. On the other hand, when $\mu$ is
sufficiently large negative, there exists a bilateral solution with the
opposite phase; i.e., $y_{\text{upper}}$ occurs for $n$ even. We refer to the
latter as an \emph{out-of-phase} bilateral solution.

Continuing the notation as in (\ref{eqn:bilateral}), the out-of-phase
bilateral solution is characterized by the following equations%
\begin{equation}
y_{\text{upper}}=B\,y_{\text{lower}}+\left(  \mu-\delta\right)
c,\,\,\,\,\,\,\,y_{\text{lower}}=A\,y_{\text{upper}}+\left(  \mu
+\delta\right)  c. \label{eqn:out_eqn}%
\end{equation}
Note the signs of the terms $\pm\delta\,c$ are different in the above
equations and (\ref{eqn:bilateral_eqn}). Solving the above equations yields%
\begin{equation}
y_{\text{upper}}\left(  \mu,\delta\right)  =x_{\text{upper}}\left(
\mu\right)  -\delta\,s_{\text{upper}},\,\,\,\,\,\,\,y_{\text{lower}}\left(
\mu,\delta\right)  =x_{\text{lower}}\left(  \mu\right)  +\delta
\,s_{\text{lower}}, \label{eqn:out_sol}%
\end{equation}
where $x_{\text{upper}}$ and $x_{\text{lower}}$ have been computed in Equation
(\ref{eqn:x_u&l}), and the $\mu$-independent shifts $s_{\text{upper}}$,
$s_{\text{lower}}$ have been computed in Equation (\ref{eqn:s_u&l}). Note that
the out-of-phase and in-phase bilateral solutions are parallel to each other
because they are both parallel to the solution $x_{\text{upper}}\left(
\mu\right)  $, $x_{\text{lower}}\left(  \mu\right)  $. One can show that the
out-of-phase bilateral solution exists only when $\mu\,X_{\text{upper}%
}^{\left(  1\right)  }\geq\delta\,s_{\text{upper}}^{\left(  1\right)
}$. As required in Condition 6, $X_{\text{upper}}^{\left(  1\right)
}<0$. Therefore, $\mu$ has to be sufficiently large negative for the
out-of-phase bilateral solution to exist.

\end{document}